\documentclass[11pt]{amsart}

\makeatletter
\renewcommand{\author}[2][]{%
  \ifx\@empty\authors
    \gdef\authors{#2}%
  \else
    \g@addto@macro\authors{\and#2}%
    \g@addto@macro\addresses{\author{}}%
  \fi
  \@ifnotempty{#1}{%
    \ifx\@empty\shortauthors
      \gdef\shortauthors{\footnotesize\scshape#1}%
    \else
      \g@addto@macro\shortauthors{\and\footnotesize\scshape#1}%
    \fi
  }}
\edef\author{\@nx\@dblarg
  \@xp\@nx\csname\string\author\endcsname}
\let\shortauthors\@empty   \let\authors\@empty
\def\maketitle{\par
  \@topnum\z@%
  \@setcopyright
  \thispagestyle{firstpage}%
  \uppercasenonmath\shorttitle
  \ifx\@empty\shortauthors \let\shortauthors\shorttitle
  \else \andify\shortauthors
  \fi
  \@maketitle@hook
  \begingroup
  \@maketitle
  \toks@\@xp{\shortauthors}\@temptokena\@xp{\shorttitle}%
  \toks4{\def\\{ \ignorespaces}}%
  \edef\@tempa{%
    \@nx\markboth{\the\toks4
      \@nx{\the\toks@}}{\the\@temptokena}}%
  \@tempa
  \endgroup
  \c@footnote\z@
  \@cleartopmattertags}
\def\@setauthors{%
  \begingroup
  \def\thanks{\protect\thanks@warning}%
  \trivlist
  \centering\@topsep30\p@\relax
  \advance\@topsep by -\baselineskip
  \item\relax
  \author@andify\authors
  \def\\{\protect\linebreak}%
  \small\scshape\authors%
  \ifx\@empty\contribs
  \else
    ,\penalty-3 \space \@setcontribs
    \@closetoccontribs
  \fi
  \endtrivlist
  \endgroup}
\renewenvironment{abstract}{%
  \ifx\maketitle\relax
    \ClassWarning{\@classname}{Abstract should precede
      \protect\maketitle\space in AMS document classes; reported}%
  \fi
  \global\setbox\abstractbox=\vtop \bgroup
    \normalfont\Small
    \list{}{\labelwidth\z@
      \leftmargin3pc \rightmargin\leftmargin
      \listparindent\normalparindent \itemindent\z@
      \parsep\z@ \@plus\p@
      }%
    \item[\hskip\labelsep\bfseries\abstractname.]%
}{%
  \endlist\egroup
  \ifx\@setabstract\relax \@setabstracta \fi}

\renewcommand{\@secnumfont}{\bfseries}
\def\section{\@startsection{section}{1}%
  \z@{.7\linespacing\@plus\linespacing}{.5\linespacing}%
  {\normalfont\bfseries\centering}}

\def\subsection{\@startsection{subsection}{2}%
  \z@{.5\linespacing\@plus.7\linespacing}{-.5em}%
  {\normalfont\itshape}}
\def\@setaddresses{\par
  \nobreak \begingroup
\footnotesize
  \def\author##1{\nobreak\addvspace\bigskipamount}%
  \def\\{\unskip, \ignorespaces}%
  \interlinepenalty\@M
  \def\address##1##2{\begingroup
    \par\addvspace\bigskipamount\indent
    \@ifnotempty{##1}{(\ignorespaces##1\unskip) }%
    {\ignorespaces##2}\par\endgroup}%
  \def\curraddr##1##2{\begingroup
    \@ifnotempty{##2}{\nobreak\indent\curraddrname
      \@ifnotempty{##1}{, \ignorespaces##1\unskip}\/:\space
      ##2\par}\endgroup}%
  \def\email##1##2{\begingroup
    \@ifnotempty{##2}{\nobreak\indent\emailaddrname
      \@ifnotempty{##1}{, \ignorespaces##1\unskip}\/:\space
      \ttfamily##2\par}\endgroup}%
  \def\urladdr##1##2{\begingroup
    \def~{\char`\~}%
    \@ifnotempty{##2}{\nobreak\indent\urladdrname
      \@ifnotempty{##1}{, \ignorespaces##1\unskip}\/:\space
      \ttfamily##2\par}\endgroup}%
  \addresses
  \endgroup}
\makeatother

\usepackage{lmodern}
\usepackage{amsmath, amsthm, amssymb, amsfonts}
\usepackage[normalem]{ulem}
\usepackage{hyperref}

\usepackage{verbatim} 
\usepackage{longtable}

\usepackage{mathtools}

\usepackage{tikz}
\usetikzlibrary{decorations.pathmorphing}
\tikzset{snake it/.style={decorate, decoration=snake}}

\usepackage{caption}

\usepackage{tikz-cd}
\usetikzlibrary{arrows}

\theoremstyle{plain}

\theoremstyle{definition}

\theoremstyle{remark}

\newcommand{\BC}{{\mathbb{C}}}

\newcommand{\BQ}{{\mathbb{Q}}}
\newcommand{\BR}{{\mathbb{R}}}

\newcommand{\BZ}{{\mathbb{Z}}}

\newcommand{\CM}{{\mathcal M}}

\DeclareFontFamily{OT1}{rsfs}{}
\DeclareFontShape{OT1}{rsfs}{n}{it}{<-> rsfs10}{}
\DeclareMathAlphabet{\curly}{OT1}{rsfs}{n}{it}

\def\secskip{\vspace{.5\linespacing plus.7\linespacing}}
\newcommand\subsectionempty{\subsection{\texorpdfstring{\!\!\!}{}}}

\begin{document}
\title[$P=W$ for Lagrangian fibrations]{$P=W$ for Lagrangian fibrations and degenerations of hyper-K\"ahler manifolds}
\date{\today}

\author{Andrew Harder}
\address{Lehigh University, Department of Mathematics}
\email{anh318@lehigh.edu}

\author{Zhiyuan Li}
\address{Fudan University, Shanghai Center for Mathematical Sciences}
\email{zhiyuan\_li@fudan.edu.cn}

\author{Junliang Shen}
\address{Massachusetts Institute of Technology, Department of Mathematics}
\email{jlshen@mit.edu}

\author{Qizheng Yin}
\address{Peking University, Beijing International Center for Mathematical Research}
\email{qizheng@math.pku.edu.cn}

\thanks{Z.~L. was supported by the NSFC grants 11731004, 11771086, and the Shu Guang Project 17GS01; Q.~Y.~was supported by the NSFC grants 11701014 and 11831013.}

\begin{abstract}
We identify the perverse filtration of a Lagrangian fibration with the monodromy weight filtration of a maximally unipotent degeneration of compact hyper-K\"ahler manifolds.
\end{abstract}

\maketitle

\subsectionempty
Throughout, we work over the complex numbers $\BC$. Let $M$ be an irreducible holomorphic symplectic variety or, equivalently, a projective hyper-K\"ahler manifold. Assume that it admits a (holomorphic) Lagrangian fibration $\pi:M\rightarrow B$. The perverse $t$-structure on the constructible derived category $D^b_c(B, \BQ)$ induces a perverse filtration on the cohomology of~$M$,
\[P_\bullet H^*(M, \BQ).\]
We refer to \cite{dCHM1, SY} for the conventions of the perverse filtration.

\subsectionempty
Let $f: \mathcal{M} \to \Delta$ be a projective degenerating family of hyper-K\"ahler manifolds over the unit disk. For $t \in \Delta^*$, let $N$ denote the logarithmic monodromy operator on $H^2(\mathcal{M}_t, \BQ)$. The degeneration $f: \mathcal{M} \to \Delta$ is called of type III if
\[N^2 \neq 0, \quad N^3 = 0.\]
By \cite[Proposition 7.14]{KLSV}, this is equivalent to having maximally unipotent monodromy. See the rest of \cite{KLSV} and also \cite{GKLR, Nag} for more discussions on degenerations of hyper-K\"ahler manifolds.

Let
\[\big(H^*_{\mathrm{lim}}(\BQ), W_\bullet H^*_{\mathrm{lim}}(\BQ), F_\bullet H^*_{\mathrm{lim}}(\BC)\big)\]
denote the limiting mixed Hodge structure\footnote{Similarly as the perverse filtration, we consider the Hodge filtration as an increasing filtration.} associated with $f: \mathcal{M} \to \Delta$.
In this short note, we prove the following result relating the perverse and the monodromy weight filtrations.

\subsection{Theorem} \label{mainthm}
--- \emph{For any Lagrangian fibration $\pi: M \to B$, there exists a type III projective degeneration of hyper-K\"ahler manifolds $f: \mathcal{M} \to \Delta$ with $\mathcal{M}_t$ deformation equivalent to $M$ for all~$t \in \Delta^*$, such that}
\begin{equation} \label{P=W}
P_kH^*(M, \BQ) = W_{2k}H^*_{\mathrm{lim}}(\BQ) = W_{2k + 1}H^*_{\mathrm{lim}}(\BQ)
\end{equation}
\emph{through an identification of the cohomology rings $H^*(M, \BQ) = H^*_{\mathrm{lim}}(\BQ)$.}

\secskip
Since $M$ and $\mathcal{M}_t$ are deformation equivalent and hence diffeomorphic, they share the same cohomology. The limiting mixed Hodge structure can be viewed as supported on the cohomology of $\mathcal{M}_t$, which provides the required identification $H^*(M, \BQ) = H^*_{\mathrm{lim}}(\BQ)$. This identification will be built into the construction of the degeneration $f: \mathcal{M} \to \Delta$.

\subsectionempty
Theorem \ref{mainthm} was previously conjectured by the first author in \cite[Conjecture 1.4]{Har} and proven in the case of $K3$ surfaces.

The interaction between the perverse and the weight filtrations for certain (noncompact) hyper-K\"ahler manifolds was first discovered by de~Cataldo, Hausel, and Migliorini \cite{dCHM1}, which is now referred to as the $P=W$ conjecture. More precisely, the~$P=W$ conjecture identifies the perverse filtration of a Hitchin fibration with the weight filtration of the mixed Hodge structure of the corresponding character variety through Simpson's nonabelian Hodge theory~\cite{Simp}. Theorem \ref{mainthm} can be viewed as a direct analogue of this conjecture.

\subsectionempty
Theorem \ref{mainthm} also offers conceptual explanations to the main results in \cite{SY}. As is remarked in~\cite[Introduction]{Har}, a recent result of Soldatenkov \cite[Theorem 3.8]{Sol} shows that limiting mixed Hodge structure for type III degenerations is of Hodge--Tate type.\footnote{This parallels the fact that the mixed Hodge structure of character varieties is of Hodge--Tate type; see~\cite{Shende}.} In particular, we have
\[\dim_\BQ \mathrm{Gr}_{2i}^WH^{i + j}_{\mathrm{lim}}(\BQ) = \dim_\BC \mathrm{Gr}_{i}^FH^{i + j}_{\mathrm{lim}}(\BC).\]
Coupled with the equalities (by \eqref{P=W} and the definition of the limiting Hodge filtration)
\begin{gather*}
\dim_\BQ \mathrm{Gr}_{i}^PH^{i + j}(M, \BQ) = \dim_\BQ \mathrm{Gr}_{2i}^WH^{i + j}_{\mathrm{lim}}(\BQ), \\
\dim_\BC \mathrm{Gr}_{i}^FH^{i + j}_{\mathrm{lim}}(\BC) = \dim_\BC \mathrm{Gr}_{i}^FH^{i + j}(\mathcal{M}_t, \BC) = \dim_\BC \mathrm{Gr}_{i}^FH^{i + j}(M, \BC),
\end{gather*}
this yields the ``Perverse $=$ Hodge'' equality in \cite[Theorem 0.2]{SY},
\begin{equation*}
\dim_\BQ \mathrm{Gr}_{i}^PH^{i + j}(M, \BQ) = \dim_\BC \mathrm{Gr}_{i}^FH^{i + j}(M, \BC).
\end{equation*}
See \cite[Section 0.4]{SY} for various applications of this equality.

Moreover, the $P=W$ identity \eqref{P=W} implies the multiplicativity of the perverse filtration 
\begin{equation*}
\cup: P_kH^d(M, \BQ) \times P_{k'}H^{d'}(M, \BQ) \to P_{k+k'}H^{d+d'}(M, \BQ)
\end{equation*}
through the general fact that the monodromy weight filtration is multiplicative. The latter follows from a combination of results of Fujisawa and Steenbrink. Fujisawa \cite[Lemma 6.16]{Fuji} proved that the wedge product on the relative logarithmic de Rham complex of a projective semistable degeneration induces a cup product on the hypercohomology groups which respects a particular weight filtration. In the much earlier work \cite[Section 4]{Ste}, Steenbrink identified the hypercohomology of the relative logarithmic de Rham complex with the cohomology of the nearby fiber, in such a way that the cup product matches the topological cup product and the weight filtration corresponds to the monodromy weight filtration. This recovers \cite[Theorem~A.1]{SY}.

As the proof of Theorem \ref{mainthm} uses the same ingredients as in \cite{SY}, the new way of deriving these results is not logically independent.

\subsectionempty
We now prove Theorem \ref{mainthm} and we make free use of the statements in \cite{SY}. To fix some notation, let~$\pi: M \to B$ be a Lagrangian fibration with $\dim M = 2\dim B = 2n$. The second cohomology group $H^2(M, \BZ)$ (resp.~$H^2(M, \BQ)$) is equipped with the Beauville--Bogomolov--Fujiki quadratic form $q_M(-)$ of signature $(3, b_2(M) - 3)$, where $b_2(M)$ is the second Betti number of~$M$.

Let $\eta \in H^2(M, \BQ)$ be a~$\pi$-relative ample class, and let $\beta \in H^2(M, \BQ)$ be the pullback of an ample class on $B$. We have $q_M(\beta) = 0$ and, by taking $\BQ$-linear combinations of $\eta$ and~$\beta$, we may assume $q_M(\eta) = 0$. Note that in this case, we have $b_2(M) \geq 4$.

\subsectionempty
Consider the following operators on the cohomology $H^*(M, \BQ)$,
\[L_\eta(-) = \eta \cup -, \quad L_\beta(-) = \beta \cup -.\]
In \cite[Section 3.1]{SY}, it was shown that $L_\eta$ and $L_\beta$ form $\mathfrak{sl}_2$-triples $(L_\eta, H_\eta, \Lambda_\eta)$ and $(L_\beta, H_\beta, \Lambda_\beta)$ which generate an $\mathfrak{sl}_2 \times \mathfrak{sl}_2$-action on $H^*(M, \BQ)$. The action induces a weight decomposition
\begin{equation} \label{split}
H^*(M, \BQ) = \bigoplus_{i, j} P^{i, j}
\end{equation}
with 
\[
H_\eta|_{P^{i, j}} = (i - n)\,\mathrm{id}, \quad H_\beta|_{P^{i, j}} = (j - n)\,\mathrm{id}.
\]

A key observation in \cite[Proposition 1.1]{SY} is that \eqref{split} provides a canonical splitting of the perverse filtration $P_\bullet H^*(M, \BQ)$. More precisely, we have
\begin{equation} \label{perv}
P_kH^d(M, \BQ) = \bigoplus_{\substack{i + j= d \\ i \leq k}}P^{i, j}.
\end{equation}

\subsectionempty
The $\mathfrak{sl}_2 \times \mathfrak{sl}_2$-action above is part of a larger Lie algebra action on $H^*(M, \BQ)$ introduced by Looijenga--Lunts \cite[Section 4]{LL} and Verbitsky \cite{Ver90, Ver96}. The Looijenga--Lunts--Verbitsky algebra
\[\mathfrak{g} \subset \mathrm{End}\big(H^*(M, \BQ)\big)\]
is defined to be the Lie subalgebra generated by all $\mathfrak{sl}_2$-triples $(L_\omega, H, \Lambda_\omega)$ with $\omega \in H^2(M, \BQ)$ such that $L_\omega(-) = \omega \cup -$ satisfies hard Lefschetz.

Given a $\BQ$-vector space $V$ equipped with a quadratic form $q$, we define the Mukai extension
\[\widetilde{V} = V \oplus \BQ^2, \quad \tilde{q} = q \oplus \left(\begin{smallmatrix}0&-1\\-1&0\end{smallmatrix}\right).\]
Looijenga--Lunts \cite[Proposition 4.5]{LL} and Verbitsky \cite[Theorem 1.4]{Ver96} showed independently
\[\mathfrak{g} \simeq \mathfrak{so}(\widetilde{H}^2(M, \BQ), \tilde{q}_M), \quad \mathfrak{g}_\BR \simeq \mathfrak{so}(4, b_2(M) - 2).\]
Here the statement with $\BQ$-coefficients is taken from \cite[Theorem 2.7]{GKLR}. Moreover, there is a weight decomposition $\mathfrak{g} = \mathfrak{g}_{-2} \oplus \mathfrak{g}_0 \oplus \mathfrak{g}_2$ with
\begin{equation} \label{giso}
\mathfrak{g}_{-2} \simeq H^2(M, \BQ), \quad \mathfrak{g}_0 \simeq \mathfrak{so}(H^2(M, \BQ), q_M) \oplus \langle H\rangle, \quad \mathfrak{g}_{2} \simeq H^2(M, \BQ).
\end{equation}

Another relevant Lie algebra is generated by the $\mathfrak{sl}_2$-triples associated with $\eta$, $\beta$, and a third element $\rho \in H^2(M, \BQ)$ satisfying
\begin{equation*} \label{orth}
q_M(\rho) > 0, \quad q_M(\eta, \rho) = q_M(\beta, \rho) = 0.
\end{equation*}
Such a $\rho$ exists by the signature $(3, b_2(M) - 3)$ of $q_M$. Let $\mathfrak{g}_\rho \subset \mathfrak{g}$ denote this Lie subalgebra and let
\[V_\rho = \langle\eta, \beta, \rho\rangle \subset H^2(M, \BQ).\]
By \cite[Corollary 2.5]{SY} complemented with the argument in \cite[Theorem 2.7]{GKLR}, we have
\begin{equation} \label{so5}
\mathfrak{g}_\rho \simeq \mathfrak{so}(\widetilde{V}_\rho,\, \tilde{q}_M|_{\widetilde{V}_\rho}).
\end{equation}
The $\mathfrak{g}_\rho$-action on $H^*(M, \BQ)$ induces the same weight decomposition as \eqref{split}; see \cite[Section 3.1]{SY}. 

\subsectionempty
Recall the natural isomorphism $\bigwedge^2 H^2(M,\BQ)\simeq \mathfrak{so}(H^2(M,\BQ), q_M)$ defined by
$$a\wedge b\mapsto \frac{1}{2}q_M(a,-)\,b-\frac{1}{2}q_M(b,-)\,a.$$
As in \cite[Lemma 4.1]{Sol}, we obtain a nilpotent operator $N_{\beta,\rho}=\beta\wedge \rho \in \mathfrak{so}(H^2(M,\BQ), q_M)$ whose action on $H^2(M, \BQ)$ satisfies
\[\mathrm{Im}(N_{\beta,\rho})=\langle\beta,\rho\rangle, \quad \mathrm{Im}(N^2_{\beta,\rho})=\langle\beta\rangle, \quad N^3_{\beta,\rho}=0.\]
By \cite[Lemma 3.9]{KSV} and the assumption $q_M(\beta,\rho)=0$, we can further identify $N_{\beta,\rho}$ with the commutator $[L_\beta,\Lambda_\rho] \in \mathfrak{g}_0$ through the isomorphisms \eqref{giso}. Note that $N_{\beta,\rho} = [L_\beta, \Lambda_\rho] \in \mathfrak{g}_\rho$.

In the two remaining sections, we show that the nilpotent operator $N_{\beta, \rho}$ induces an~$\mathfrak{sl}_2$-triple whose weight decomposition splits both the perverse filtration $P_\bullet H^*(M, \BQ)$ and the monodromy weight filtration of a degeneration $f: \CM \to \Delta$. This completes the proof of Theorem \ref{mainthm}.

\subsectionempty
The construction of a degeneration $f:\mathcal{M}\rightarrow \Delta$ with logarithmic monodromy $N_{\beta,\rho}$ is precisely \cite[Theorem 4.6]{Sol}. While the original statement requires $b_2(M)\geq 5$ to ensure the existence of an element~$\beta \in H^2(M, \BQ)$ with $q_M(\beta) = 0$, in our situation $\beta$ is readily given by the Lagrangian fibration~$\pi: M \to B$. From the proof of \cite[Theorem~4.6]{Sol}, it suffices to find an element~$h\in H^2(M,\mathbb{Z})$ satisfying 
\[
q_M(h)>0,\quad q_M(\beta, h)=q_M(\rho, h)=0
\]
in order to obtain nilpotent orbits $(N_{\beta,\rho}, x)$ with $x \in \widehat{\mathcal{D}}_h$ as in \cite[Definition~4.3]{Sol}.\footnote{Here $\widehat{\mathcal{D}}_h$ is the extended polarized period domain with respect to $h \in H^2(M, \BZ)$.} These nilpotent orbits eventually provide the required degeneration $f: \mathcal{M} \to \Delta$ through global~Torelli. Now since $q_M$ is of signature $(3, b_2(M) - 3)$ while $q_M|_{V_\rho}$ is only of signature $(2, 1)$ (recall that~$b_2(M) \geq 4$), such an $h$ exists.

By Jacobson--Morozov, the nilpotent operator $N_{\beta, \rho} \in \mathfrak{g}_\rho$ is part of an $\mathfrak{sl}_2$-triple which we denote $(L_N = N_{\beta, \rho}, H_N, \Lambda_N)$. Consider the action of this $\mathfrak{sl}_2$ on $H^*(M, \BQ)$ and the associated weight decomposition
\begin{equation} \label{wsplit}
H^*(M, \BQ) = \bigoplus_{d, m}W_m^d
\end{equation}
with $H_N|_{W_m^d} = m\,\mathrm{id}$. By the definition of the monodromy weight filtration, we have
\begin{equation} \label{weight}
W_kH^d_{\mathrm{lim}}(\BQ) = \bigoplus_{d - m \leq k} W_m^d.
\end{equation}

\subsectionempty 
Finally, we match the perverse decomposition \eqref{split} with the weight decomposition \eqref{wsplit}. As both decompositions are defined over $\BQ$, it suffices to work with $\BC$-coefficients.

We recall some basic facts about $\mathfrak{so}(5,\BC)$-representations. Let $V$ be a $\BC$-vector space admitting three $\mathfrak{sl}_2$-actions $(L_1, H, \Lambda_1)$, $(L_2, H, \Lambda_2)$, and $(L_3, H, \Lambda_3)$ which generate an $\mathfrak{so}(5,\BC)$-action. More concretely, the operators
\[
L_s,\, \Lambda_s,\, K_{st}=[L_s, \Lambda_t],\, H, \quad \text{for } s,t \in \{1,2,3\}
\]
satisfy the relations (2.1) in \cite{Ver90}. We consider the Cartan subalgebra
\[
\mathfrak{h}= \langle H, -\sqrt{-1}K_{23} \rangle \subset \mathfrak{so}(5, \BC)
\]
and the associated weight decomposition 
\[
V= \bigoplus_{i,j} V^{i,j}
\]
with
\[
H|_{V^{i,j}} = (i+j -2n)\,\mathrm{id}, \quad (-\sqrt{-1}K_{23})|_{V^{i,j}} = (i-j)\,\mathrm{id}.
\]

We define a nilpotent operator
\begin{equation*}
L_N = \left[\frac{1}{2}L_2- \frac{\sqrt{-1}}{2}L_3,\, \Lambda_1 \right] = -\frac{1}{2} K_{12}+\frac{\sqrt{-1}}{2}K_{13} \in \mathfrak{so}(5, \BC),
\end{equation*}
which induces an $\mathfrak{sl}_2$-triple $(L_N, H_N, \Lambda_N)$ with
\[
\Lambda_N = \left[-\frac{1}{2}L_2- \frac{\sqrt{-1}}{2}L_3,\, \Lambda_1 \right] = \frac{1}{2} K_{12}+\frac{\sqrt{-1}}{2}K_{13}, \quad H_N = \sqrt{-1}K_{23}.
\]
In particular, we have $H_N|_{V^{i,j}} = (j-i)\,\mathrm{id}$. The weight decomposition with respect to this~$\mathfrak{sl}_2$-action then takes the form
\[
V= \bigoplus_{m} V_m^d, \quad V_m^d = \bigoplus_{\substack{i + j = d \\ j - i = m}} V^{i,j}
\]
with $H_N|_{V_m^d} = m\,\mathrm{id}$.

In our geometric situation, let $V$ be the total cohomology $H^*(M, \BC)$. We consider the three operators $L_1, L_2, L_3$ determined by
\[
L_1 = L_\rho,\quad \frac{1}{2}L_2+ \frac{\sqrt{-1}}{2}L_3 =L_\eta,\quad \frac{1}{2}L_2-\frac{\sqrt{-1}}{2}L_3 = L_\beta
\]
which induce a representation of $\mathfrak{so}(5, \BC)$ by \eqref{so5}. In particular, we have $V^{i,j} = P^{i,j}_\BC$. Moreover, the nilpotent operator $L_N$ is exactly $N_{\beta, \rho} = [L_\beta, \Lambda_\rho]$. We conclude from \eqref{perv} and \eqref{weight}
\[
P_kH^d(M, \BC) = \bigoplus_{\substack{i + j= d \\ i \leq k}} P^{i,j}_\BC = \bigoplus_{\substack{i + j = d \\ j - i = m \\ d - m \leq 2k}} V^{i,j} = \bigoplus_{d - m \leq 2k} V_m^d = \bigoplus_{d - m \leq 2k} W_{m, \BC}^d = W_{2k}H^d_{\mathrm{lim}}(\BC).
\]

\end{document}